\documentclass[letterpaper, 10 pt, conference]{ieeeconf}  % Comment this line out
\usepackage{times}
\usepackage{algorithm, algorithmic}
\usepackage{amsmath}
\usepackage{amssymb}
\usepackage{array}
\usepackage[british]{babel}
\usepackage{booktabs}
\usepackage{dsfont}
\usepackage{mathtools}
\usepackage{multirow}
\usepackage{xcolor}
\usepackage{cuted}
\usepackage{accents}

\newtheorem{thm}{Theorem}[section]
\newtheorem{lem}[thm]{Lemma}

\newcommand{\xalgm}{\mathcal{X}_{\rm Alg}}

\newcommand{\tqn}{\mathcal{T}_\mathcal{Q}^{(N)}}
\newcommand{\qn}{\mathcal{Q}^{(N)}}
\newcommand{\oqn}{{\qn}^\star}
\newcommand{\barqn}{\underbar{\text{$\mathcal{Q}$}}^{(N)}}
\DeclareMathOperator*{\argmin}{arg\,min}
\DeclareMathOperator*{\argmax}{arg\,max}
\mathchardef\mhyphen="2D

\IEEEoverridecommandlockouts                              % This command is only
\title{\LARGE \bf
Learning Q-function approximations for hybrid control problems
}

\author{Sandeep Menta, Joseph Warrington, John Lygeros and Manfred Morari% <-this % stops a space
\thanks{S.Menta, J.Warrington and J.Lygeros are with Automatic Control Lab, ETH Z\"urich, Physikstrasse 3, 8092 Z\"urich, Switzerland {\tt\small email: smenta@ethz.ch; joe.warrington@gmail.com; jlygeros@ethz.ch}}%
\thanks{M.Morari is with Elec. and Systems Engineering, Univ. of Pennsylvania, 220 S. 33rd St, Philadelphia, PA 19104, United States {\tt\small email: morari@seas.upenn.edu}}%
}

\begin{document}

\maketitle
\thispagestyle{empty}
\pagestyle{empty}

%%%%%%%%%%%%%%%%%%%%%%%%%%%%%%%%%%%%%%%%%%%%%%%%%%%%%%%%%%%%%%%%%%%%%%%%%%%%%%%%
\begin{abstract}
The main challenge in controlling hybrid systems arises from having to consider an exponential number of sequences of future modes to make good long-term decisions. Model predictive control (MPC) computes a control action through a finite-horizon optimisation problem. A key ingredient in this problem is a terminal cost, to account for the system's evolution beyond the chosen horizon. A good terminal cost can reduce the horizon length required for good control action and is often tuned empirically by observing performance. We build on the idea of using $N$-step $Q$-functions $(\qn)$ in the MPC objective to avoid having to choose a terminal cost. We present a formulation incorporating the system dynamics and constraints to approximate the optimal $\qn$-function and algorithms to train the approximation parameters through an exploration of the state space. We test the control policy derived from the trained approximations on two benchmark problems through simulations and observe that our algorithms are able to learn good $\qn$-approximations for high dimensional hybrid systems based on a relatively small data-set. Finally, we compare our controller's performance against that of Hybrid MPC in terms of computation time and closed-loop cost.
\end{abstract}

%%%%%%%%%%%%%%%%%%%%%%%%%%%%%%%%%%%%%%%%%%%%%%%%%%%%%%%%%%%%%%%%%%%%%%%%%%%%%%%%
\section{Introduction}
Hybrid systems have been used to describe a range of real-world systems like solar energy plants \cite{Solar}, bipedal robots \cite{bipedal} and vehicle dynamics \cite{traction_control}. Hybrid systems allow one to use different dynamics at different parts of the state space at the price of having to include integer states or inputs complicating hybrid control problems. A popular control approach is hybrid Model Predictive Control (HMPC) that optimises over $N$-step trajectories with a terminal cost and/or constraints to account for the system evolution beyond the planning horizon. A well chosen terminal cost or constraint can reduce $N$, hence the computational effort required, but no general approach for making this choice exists.  \textit{Explicit MPC}, \cite{explicit_mpc} methods reduce the online optimisation effort by computing the explicit feedback law defined over polyhedral partitions of states offline. For larger systems and longer horizon lengths the number of partitions can become prohibitive for online use. To reduce the computation Neural network (NN) based methods have emerged to approximate the explicit MPC laws by incorporating the model into their training procedures \cite{Chen_neural} or using specialised layers \cite{maddalena2019neural}.

We follow the approach presented in \cite{menta2020} to use an $N$-step extension of the $Q$-function, $\qn$, that takes as its arguments the first $N$ states and actions. The optimal control policy can then be computed by minimising the optimal $N$-step $Q$-function $(\oqn)$, a process comparable to hybrid MPC. Our approach is to learn approximations of $\oqn$ as a point-wise maximum of several lower bounding functions to $\oqn$. The lower bounding functions are built by incorporating system dynamics and constraints. Our learning algorithm iteratively improves the approximation by adding tighter lower bounding functions at states chosen based on a metric of improvement. The advantage of our formulation is that it allows us to extract a control policy by solving a mixed-integer quadratic program (MIQP) of complexity comparable to HMPC instead of the computationally expensive \emph{quadratically constrained} MIQCQP which results from the formulation proposed in \cite{menta2020}.

We also test the necessity of utilising model information to achieve good control by using off-the-shelf NN training techniques and architectures, which do not use any model information, to learn the optimal control actions using the current state information. We compare the performance of our methods with HMPC and explicit controllers based on off-the-shelf NN training methods. Our simulation results on benchmark problems suggest that off-the-shelf NN training methods are not a good choice for hybrid control tasks and the controller derived from our approximations performs comparably to HMPC in terms of computational ease as well as the quality of the computed control actions.

The $\qn$-function and the control policy it generates are discussed in Section \ref{sec:NStepQ}. The algorithm used to approximate the $\oqn$ is described in Section \ref{sec:NovelFormulation} and the simulation results on benchmark examples are presented in Section \ref{sec:NumEx}.
%%%%%%%%%%%%%%%%%%%%%%%%%%%%%%%%%%%%%%%%%%%%%%%%%%%%%%%%%%%%%%%%%%%%%%%%%%%%%%%%
\section{N-step $\mathcal{Q}$-Function}\label{sec:NStepQ}
\subsection{Optimal $V$ function}
We consider the problem of optimal control for hybrid systems represented as time-invariant mixed logical dynamical (MLD) systems \cite{bemporad1999control} without binary states or control inputs. Our starting point is an infinite horizon optimal control problem whose value function $V^\star$ is given by:
\allowdisplaybreaks{
\begin{subequations} \label{eq:Vdef}
	\begin{align}
	V^\star(x) := & \nonumber\\
	\min_{\substack{\{x\}_{0}^\infty, \{u\}_{0}^\infty,\\ \\\{\delta\}_{0}^\infty, \{z\}_{0}^\infty}} \, &\sum_{t=0}^{\infty} \frac{\gamma^t}{2} \left( \begin{aligned}(x_t -& x_g)^\top Q (x_t - x_g) \,\, +\\ &(u_t - u_g)^\top R (u_t - u_g) \end{aligned} \right) \\
	\text{s.~t.}\quad & \begin{aligned}
		x_{t+1} = Ax_t + B_1 u_t + B_2 \delta_t + B_3 z_t \, , \,\,\,\\ t = 0, 1, \ldots \, , 
	\end{aligned}\label{eq:VdefMLDDyn}\\
	& \begin{aligned}
		E_2 \delta_t + E_3 z_t \leq E_4 x_t + E_1 u_t + E_5 \, , \,\,\,\\ t = 0, 1, \ldots \, , 
	\end{aligned}\label{eq:VdefMLDConstr} \\
	& \delta_t \in \{0,1\}^{n_\delta} \, , \,\,\, t = 0, 1, \ldots \, , \\
	& x_0 = x \, ,  \label{eq:VdefInitStateConstr}
	\end{align}
\end{subequations}
}where $x_t, x_g \in \mathbb{R}^n$ are the state and the goal state, $u_t, u_g \in \mathbb{R}^m$ are the input and the desired input, $z_t \in \mathbb{R}^{n_z}$ and $\delta_t$ are auxiliary continuous and binary variables required in the MLD framework. The shorthands $\{x\}_{0}^{N}$, are used to denote the sequences $x_0, x_1,\dots,x_N$, etc. We assume $Q \succ 0$ and $R \succ 0$. Constraints \eqref{eq:VdefMLDDyn}-\eqref{eq:VdefInitStateConstr} ensure that the system evolves along feasible trajectories, including mode switches parameterized by $\delta_t$. We assume the MLD system is ``well-posed'' \cite{bemporad1999control}, in the sense that $x_t$ and $u_t$ uniquely determine $\delta_t$ and $z_t$. The MLD constraints \eqref{eq:VdefMLDConstr} comprise the constraints that ensure the MLD system is ``well-posed" and the state and input constraints which depend only on $x_t$ and $u_t$, with $n_c$ constraints per time step. The discount factor $\gamma \in (0, 1]$  is used here to include a wider class of problems, for example, problems with periodic solutions.

\subsection{Hybrid MPC}
HMPC involves computing the optimal control actions by minimising an approximation of $V^\star(x)$. The objective of HMPC uses the first $N$-terms from \eqref{eq:Vdef} and approximates the rest of the summation via a terminal cost $V(x_N)$. This formulation using the truncated objective can be found in \cite[\S 17.4]{borrelli2017predictive}. The ideal terminal cost would be the optimal cost function $V^\star(x_N)$ itself but $V^\star (x)$ is in general impossible to compute. To reduce the effect of truncating the infinite horizon problem and using a sub-optimal terminal cost the HMPC controller is used in a receding horizon fashion.

\subsection{Optimal $N$-step $\mathcal{Q}$-function}
The notion of a $N$-step $\mathcal{Q}$-function introduced in \cite{menta2020} takes the initial state $x_0$ and a $N$-step feasible sequence of inputs and auxiliary variables $(\{u,\delta,z\}_0^N)$ as its parameters. Though redundant due to well-posedness, without loss of generality one can also include $\{x\}_1^N$ as arguments. The optimal $N$-step $\mathcal{Q}$-function $(\oqn)$ associated with this parametrization of the $\qn$-function can be obtained in the same way the optimal value function is computed:

\begin{align}
	\oqn (\{x&\}_0^{N \mhyphen 1}, \{u\}_0^{N \mhyphen 1}, \{\delta\}_0^{N \mhyphen 1}, \{z\}_0^{N \mhyphen 1}) := \nonumber\\ 
	\min_{\substack{\{x\}_{N}^\infty, \{u\}_{N}^\infty,\\ \{\delta\}_{N}^\infty, \{z\}_{N}^\infty}}& \sum_{t=0}^{\infty} \frac{\gamma^t}{2} \left( \begin{aligned}(x_t -& x_g)^\top Q (x_t - x_g) \,\, +\\ &(u_t - u_g)^\top R (u_t - u_g) \end{aligned} \right) \label{eq:NstepQdef}\\ 
	\text{s.t.~} & \text{\eqref{eq:VdefMLDDyn}-\eqref{eq:VdefInitStateConstr}} \nonumber
\end{align}
which implies that
\begin{multline}
\oqn(\{x\}_0^{N \mhyphen 1}, \{u\}_0^{N \mhyphen 1}, \{\delta\}_0^{N \mhyphen 1}, \{z\}_0^{N \mhyphen 1}) =\\ \textstyle \sum_{t=0}^{N \mhyphen 1} \frac{\gamma^t}{2} (x_t - x_g)^\top Q (x_t - x_g) +\\ \textstyle \sum_{t=0}^{N \mhyphen 1} \frac{\gamma^t}{2} (u_t - u_g)^\top R (u_t - u_g) + \gamma^N V^\star(x_N) ,\label{eq:QNStarVStar}
\end{multline}
where $\{x\}_{0}^N, \{u\}_{0}^N, \{\delta\}_{0}^N, \{z\}_{0}^{N}$ satisfy the constraints \eqref{eq:VdefMLDDyn}-\eqref{eq:VdefInitStateConstr}. To simplify notation we introduce $s = (x, u,\delta,z) \in \mathbb{R}^{n} \times \mathbb{R}^{m} \times \{0,1\}^{n_\delta} \times \mathbb{R}^{n_z}$. The state $x$ along with the inputs and auxiliary variables compatible with it can then be described by the set $\mathcal{S}(x)\coloneqq \{s = (x, u,\delta,z) \in \mathbb{R}^{n} \times \mathbb{R}^{m} \times \{0,1\}^{n_\delta} \times \mathbb{R}^{n_z} \,\,:\,\, E_2 \delta + E_3 z \leq E_4 x + E_1 u + E_5\}$ and the set of $N$-step feasible trajectories from the initial state $x$,
\begin{equation}\label{eq:DefSN}
	\mathcal{S}^{(N)}(x) \coloneqq \left\{\begin{multlined}
		\{s\}_0^{N\mhyphen1}=\{x, u, \delta, z\}_0^{N \mhyphen 1} : x_0 = x,\\
		x_{k+1} = Ax_k + B_1 u_k + B_2 \delta_k + B_3 z_k,\\
		s_k \in \mathcal{S}(x_k) \text{ for } k=0,\dots,N-1
	\end{multlined}
	\right\}
\end{equation}

Following \cite{menta2020}, a greedy control policy can be extracted from a $\qn$-function as follows:
\begin{equation}\label{eq:NstepPi}
\pi(x; \qn) \in \left[ \underset{\{s\}_0^{N \mhyphen 1} \in \mathcal{S}^{(N)}(x)}{\arg \min} \, \qn(\{s\}_0^{N \mhyphen 1})\right]_{u_0}.
\end{equation}
The notation $[\cdots]_{u_0}$ shows that the control input $u_0$ is taken from the optimal $s_0$, assuming an optimum is attained. We can conclude from \eqref{eq:NstepQdef} that $\pi(x; {Q^{(N)}}^\star)$ is an optimal policy that gives the optimal infinite-horizon ``closed-loop'' cost $V^\star(x)$ when applied recursively.
%
%\textcolor{red}{By the definition of $V^\star$ in  \eqref{eq:Vdef}, any other $Q^{(N)}$ will generally induce closed-loop costs greater than $V^\star(x)$.}
%
\section{Novel approximate $\oqn$-formulation}\label{sec:NovelFormulation}
The ease of extracting the control policy in \eqref{eq:NstepPi} depends heavily on the structure of the $\qn$-function. Here we construct an approximation of $\oqn$ as a point-wise maximum of several lower bounding functions referred to as ``cuts". Our formulation of the lower bounding functions builds on \cite{menta2020} and addresses some key difficulties arising there.

\subsection{Bellman operator for $N$-step $\qn$-functions}
The Bellman operator $\tqn$ for a generic function $\qn$ for a feasible trajectory $\{s\}_0^{N \mhyphen 1}$ starting at $x_0$ is defined as:
\begin{multline}\label{eq:tqn_operator}
\tqn \qn(\{s\}_0^{N \mhyphen 1}) = \tfrac{1}{2}(x_0 - x_g)^\top Q (x_0 - x_g) +\\ \tfrac{1}{2}(u_0 - u_g)^\top R (u_0 - u_g)  + \inf_{
	\substack{
		s_N \in \mathcal{S}(x_{N})
	}
} \gamma \qn(\{s\}_1^N) \,. 
\end{multline}
From \eqref{eq:NstepQdef} and \eqref{eq:tqn_operator}, for all $x_0 \in \mathbb{R}^n$, $\{s\}_0^{N \mhyphen 1} \in \mathcal{S}^{(N)}(x_0)$ we have $\oqn$ is a fixed point of $\tqn$, that is
$\tqn\oqn(\{s\}_0^{N \mhyphen 1})=\oqn(\{s\}_0^{N \mhyphen 1}).$
Moreover, $\tqn$ is monotonic, that is, if $\qn_a(\{s\}_0^{N \mhyphen 1}) \leq \qn_b(\{s\}_0^{N \mhyphen 1})$ for all $x_0 \in \mathbb{R}^n$, $\{s\}_0^{N \mhyphen 1} \in \mathcal{S}^{(N)}(x_0)$, then
\begin{multline*}
\tqn \qn_a(\{s\}_0^{N \mhyphen 1}) \leq \tqn \qn_b(\{s\}_0^{N \mhyphen 1}),\\ \forall x_0 \in \mathbb{R}^n, \{s\}_0^{N \mhyphen 1} \in \mathcal{S}^{(N)}(x_0).
\end{multline*}
\subsection{Reformulated Benders' cuts}
The approximation of $\oqn$ constructed as a point-wise maximum of $I+1$ cuts $q_0,q_1,\dots,q_I$ is given by:
\begin{equation}\label{eq:Qpw_max}
	\qn_I(\{s\}_0^{N \mhyphen 1}) \coloneqq \max_{i=0,\dots,I}\left\{q_i(\{s\}_0^{N \mhyphen 1})\right\},
\end{equation}
where each of the cuts satisfies,
\begin{multline}
	q_i(\{s\}_0^{N \mhyphen 1}) \leq \oqn(\{s\}_0^{N \mhyphen 1}),\\ \forall x_0 \in \mathbb{R}^n, \{s\}_0^{N \mhyphen 1} \in \mathcal{S}^{(N)}(x_0)
\end{multline}
This results in an approximation which lower bounds $\oqn$ everywhere. Following \cite{menta2020}, \cite{warrington2019learning} we choose the cuts to be linearly separable in their arguments:
\begin{multline}\label{eq:Separableqi}
	q_i(\{s\}_0^{N \mhyphen 1}) =  \textstyle \sum_{j=0}^{N \mhyphen 1} q_i^{x_j}(x_j) + \sum_{j=0}^{N \mhyphen 1} q_i^{u_j}(u_{j}) + \\ \textstyle \sum_{j=0}^{N \mhyphen 1} q_i^{\delta_j}(\delta_{j}) + \sum_{j=0}^{N \mhyphen 1} q_i^{z_j}(z_{j}) + c_i,
\end{multline} where $\{s\}_0^{N \mhyphen 1} \in \mathcal{S}^{(N)}(x_0)$ and each of the $q_i^{(\cdot)}(\cdot)$ terms are linear and convex quadratic functions, and $c_i$ is a constant term. Cut $q_0$ is the starting point for building the approximation:
\begin{multline}\label{eq:functionq0}
	q_0(\{s\}_0^{N \mhyphen 1}) = \textstyle \sum_{t=0}^{N \mhyphen 1} \tfrac{\gamma^{t}}{2}(x_t - x_g)^\top Q (x_t - x_g) +\\ \textstyle \sum_{t=0}^{N \mhyphen 1} \tfrac{\gamma^{t}}{2} (u_t - u_g)^\top R (u_t - u_g).
\end{multline}
From \eqref{eq:QNStarVStar}, $q_0$ is clearly a lower bounding function to $\oqn$, as $V^\star(x_N) \geq 0$.
The first step in this process is to apply the Bellman operator to the existing approximation $\qn_I$ at $\{s\}_0^{N-1} \in \mathcal{S}^{(N)}(x_0)$:
\begin{subequations}\label{eq:tqn_primal}
	\begin{align}
	&\tqn \qn_I(\{s\}_0^{N \mhyphen 1}) =\nonumber\\
	&\inf_{x_N, s_N, \alpha} \quad \left(\begin{aligned}\tfrac{1}{2}(&x_0 - x_g)^\top Q (x_0 - x_g) \,\,+\\ &\tfrac{1}{2}(u_0 - u_g)^\top R (u_0 - u_g)  + \gamma \alpha \end{aligned}\right)\\
	&\text{s.~t.}\quad x_{N} = Ax_{N\mhyphen1} + B_1 u_{N\mhyphen1} + B_2 \delta_{N\mhyphen1} + B_3 z_{N\mhyphen1}\label{eq:tqn_dyn}\\
	\label{eq:tqn_htil}&\qquad E_2 \delta_N + E_3 z_N \leq  E_1 u_N + E_4 x_N + E_5,\\
	\label{eq:tqn_qi}&\qquad q_i(\{s\}_1^{N}) \leq \alpha,\quad i=0,\dots,I
	\end{align}
\end{subequations}
where the value of $\qn_I( \{s\}_1^{N})$ is modelled by the epigraph variable $\alpha$. The form of the terms in the cuts \eqref{eq:Separableqi} comes from the objective function of the dual of \eqref{eq:tqn_primal}. We derive the dual by assigning multipliers  $\nu \in \mathbb{R}^{n}$ to \eqref{eq:tqn_dyn}, $\mu \in \mathbb{R}^{n_c}$ to \eqref{eq:tqn_htil} and $\lambda \in \mathbb{R}^{I+1}$ to \eqref{eq:tqn_qi}. The dual is then written as the operator $\mathcal{D}^{(N)}$ acting on $\qn_I$:
\begin{subequations}\label{eq:tqn_dual}
	\begin{align}
	\mathcal{D}^{(N)}&\qn_I(\{s\}_0^{N \mhyphen 1}) \coloneqq  \nonumber\\
	\sup_{\nu, \mu, \lambda} \,\,\,\, &\left\{ \begin{aligned}
	&\tfrac{1}{2}(x_0 - x_g)^\top Q (x_0 - x_g) + \textstyle \sum_{i=1}^I \lambda_{i}c_i +\\
	&\tfrac{1}{2}(u_0 - u_g)^\top R (u_0 - u_g) - {\mu}^\top E_5 + \xi(\nu, \mu, \lambda) + \\
	%& \textstyle \sum_{t=0}^{N-2}\nu_t^\top \left(
	%\begin{aligned}
	%	A x_t \,+ \,B_1 u_t \,+ &\,B_2 \delta_t + \\ 
	%	&B_3 z_t - x_{t+1}
	%\end{aligned}\right) + \\
	& \nu^\top \left(
	\begin{aligned}
		A x_{N \mhyphen 1} + B_1 u_{N \mhyphen 1} + B_2 \delta_{N \mhyphen 1} + B_3 z_{N \mhyphen 1}
	\end{aligned}\right) +  \\
	& \textstyle \sum_{i=0}^I \lambda_{i} \sum_{t=1}^{N \mhyphen 1} \left(
	\begin{aligned}
		q_i^{x_{t-1}}(x_t) &+ q_i^{u_{t-1}}(u_t) \, +\\
		q_i^{\delta_{t-1}} (\delta_t) &+ q_i^{z_{t-1}} (z_t)
	\end{aligned}\right)
	\end{aligned}\right\}\\
	\text{s.~t.} \quad & \mu \geq 0, \quad \lambda \geq 0, \quad \mathds{1}^\top \lambda = \gamma,
	\end{align}
\end{subequations}
where,
\begin{align*}
	\xi(&\nu, \mu, \lambda) \coloneqq \nonumber\\
	\inf_{s^N}&\,\,\,
	\left\{\begin{aligned}
		&\qquad \mu^\top (E_2 \delta_N + E_3 z_N - E_1 u_N - E_4 x_N) \,\,+ \nonumber \\
		& \textstyle\sum_{i=0}^I \lambda_{i} \left(
		\begin{aligned}
		&q_i^{x_{N \mhyphen 1}}(x_N) + q_i^{u_{N \mhyphen 1}}(u_N) +\\
		&q_i^{\delta_{N \mhyphen 1}}(\delta_N) + q_i^{z_{N \mhyphen 1}}(z_N) 
		\end{aligned}
		\right) - \nu^\top x_N
	\end{aligned}\right\}.
\end{align*}
The form of the new cut $q_{I+1}$ constructed using the optimal dual multipliers $(\nu^\star, \mu^\star, \lambda^\star)$ is shown in Lemma \ref{lem:new_lb_N}.
\allowdisplaybreaks{
\begin{lem}\label{lem:new_lb_N}
	The function
	\begin{multline}\label{eq:newlb_defn}
	\textstyle q_{I+1}(\{s\}_0^{N \mhyphen 1}) :=
	\sum_{t=0}^{N \mhyphen 1} q_{I+1}^{x_t}(x_t) + \sum_{t=0}^{N \mhyphen 1} q_{I+1}^{u_t}(u_{t})\,\, +\\ \textstyle \sum_{t=0}^{N \mhyphen 1} q_{I+1}^{\delta_t}(\delta_{t}) + \sum_{t=0}^{N \mhyphen 1} q_{I+1}^{z_t}(z_{t}) + c_{I+1} \, ,
	\end{multline}
	satisfies the global lower bounding property
	\begin{multline*}
	q_{I+1}(\{s\}_0^{N \mhyphen 1}) \leq \oqn (\{s\}_0^{N \mhyphen 1}) \, , \\ \forall x_0 \in \mathbb{R}^n,\,\,\{s\}_0^{N \mhyphen 1} \in \mathcal{S}^{(N)}(x_0)
	\end{multline*}
	where,
	\begin{align*}
	&q_{I+1}^{x_t}(x_t) = \tfrac{\gamma^t}{2}(x_t - x_g)^\top Q (x_t - x_g), \qquad \quad q_{I+1}^{\delta_t}(\delta_t) = 0,\\
	&q_{I+1}^{u_t}(u_t) = \tfrac{\gamma^t}{2}(u_t - u_g)^\top R (u_t - u_g), \qquad \quad q_{I+1}^{z_t}(z_t) = 0,\\
	&q_{I+1}^{N\mhyphen1}(x_{N\mhyphen1}) = \gamma q_{I+1}^{x_{N\mhyphen2}}(x_{N\mhyphen1}) + {\nu^\star}^\top A x_{N\mhyphen1} \\
	&q_{I+1}^{u_{N\mhyphen1}}(u_{N\mhyphen1}) = \gamma q_{I+1}^{u_{N\mhyphen2}}(u_{N\mhyphen1}) + {\nu^\star}^\top B_1 u_{N\mhyphen1},\\
	&q_{I+1}^{\delta_{N\mhyphen1}}(\delta_{N\mhyphen1}) = {\nu^\star}^\top B_2 \delta_{N\mhyphen1}, \qquad q_{I+1}^{z_{N\mhyphen1}}(z_{N\mhyphen1}) = {\nu^\star}^\top B_3 z_{N\mhyphen1},\\
	&c_{I+1} = \textstyle  \sum_{i=1}^{I}\lambda_{i}^\star c_i - {{\mu}^\star}^\top E_5 + \xi(\nu^\star, {\mu}^\star, \lambda^\star),
	\end{align*}
	for $t = 0,\dots,N - 2$ and the triplet $(\nu^\star, \mu^\star, \lambda^\star)$ solves problem \eqref{eq:tqn_dual} for the parameter $\{\hat{s}\}_0^{N \mhyphen 1} \in \mathcal{S}^{(N)}(\hat{x}_0)$.
\end{lem}
}

The proof is a simple adaptation of the proof of Lemma III.1 in \cite{menta2020}. Lemmas III.2-4 of \cite{warrington2019learning} carry over to the present setting but are omitted in the interest of space. We highlight that the addition of the new cut $q_{I+1}$ at the parameter $\{s\}_0^{N \mhyphen 1} \in \mathcal{S}^{(N)}(x_0)$ leads to an improvement of $\mathcal{D}^{(N)}\qn_I(\{s\}_0^{N \mhyphen 1}) - \qn_I(\{s\}_0^{N \mhyphen 1})$ in the approximation at $\{s\}_0^{N \mhyphen 1}$. Lemma \ref{lem:new_lb_N} shows that each cut is the sum of $N$ stage costs along with a penalty on the $N^{th}$ state and a constant. Thus, the generated approximation is the standard MPC objective with a terminal cost approximation that naturally incorporates the system dynamics and constraints.

\subsection{Training algorithm}\label{ss:train_algs}
\begin{algorithm}[t]
	\caption{Modified Benders algorithm for MLD systems} 
	\label{alg:QBenders_xudz}
	\begin{algorithmic}[1]
		\STATE \textbf{Inputs:} System model; training points $\xalgm := \{x^1, \ldots, x^M\}$, $I_\text{max}$, $\eta_{min}$ and $q_0(\{s\}_0^{N \mhyphen 1})$
		\STATE Set $I = 0$, $\eta^\star = \infty$
		\WHILE{$I \leq I_\text{max}$ and $\eta^\star \geq \eta_{min}$}
			\STATE $\qn_I(\cdot) \gets \max_{i=0,\ldots,I} q_i(\cdot)$
			\FOR{\textbf{each} $x \in \xalgm$} \label{algl:BellGapStartxu}
				\STATE $\{s\}_0^{N \mhyphen 1} \gets \argmin\limits_{\{s\}_0^{N \mhyphen 1} \in \mathcal{S}^{(N)}(x)} \left\{\qn_I(\{s\}_0^{N \mhyphen 1})\right\}$\label{algl:VarBum_xu}
				\STATE $\eta(\{s\}_0^{N \mhyphen 1}; \qn_I) \gets \tilde{\beta}(\{s\}_0^{N \mhyphen 1}; \qn_I)$
			\ENDFOR \label{algl:BellGapEndxu}
		%		\IF{$\max_{m=1,\ldots,M} \{\varepsilon(x_m, u_m; Q_I)\} \leq \varepsilon_\text{tol}$} \label{algl:ConvCheck2_xu}
		%		\STATE \textbf{break}
		%		\ENDIF \label{algl:ConvCheck2_2_xu}
			\STATE $x^\star,\{s^\star\}_0^{N \mhyphen 1} \gets \argmax\limits_{\substack{x \in \xalgm,\\ \{s\}_0^{N \mhyphen 1}\in \mathcal{S}^{(N)}(x)}} \left\{\eta(\{s\}_0^{N \mhyphen 1}; \qn_I)\right\}$ \label{algl:Pickx_xu}
			\IF{$\eta^\star = \eta(x^\star,\{s^\star\}_0^{N \mhyphen 1}; \qn_I) \geq \eta_{min}$}
				\STATE Add $q_{I+1}(\cdot)$ parameterized by the $(\nu^\star, {\mu}^\star, \lambda^\star)$ optimal for problem \eqref{eq:tqn_dual} with the parameter $(x^\star,\{s^\star\}_0^{N \mhyphen 1})$ to the set of $q_i(\cdot)$ functions
				\label{algl:AddLB2_xu}
			\ENDIF
			\STATE $I \gets I+1$
		\ENDWHILE
		\STATE \textbf{Output:} $\qn_I(\cdot) = \max_{i=0,\ldots,I} q_i(\cdot)$
	\end{algorithmic}
\end{algorithm}

The approximation of $\oqn$ is learnt iteratively by searching over a given set of data points $\xalgm$. At each point in $\xalgm$ the $N$-step greedy policy is evaluated using the current approximation $\qn_I$, and an improvement metric is evaluated and stored. A new cut is then added at the point which promises the highest improvement, as described in Algorithm \ref{alg:QBenders_xudz}. We propose the metric $\tilde{\beta}$ building on $\beta$ \cite{menta2020} :
\begin{align}
	&\begin{aligned}
		\beta(x;\qn_I) \coloneqq
		\mathcal{D}^{(N)} \qn_I (\{s\}_0^{N \mhyphen 1}) - \qn_I (\{s\}_0^{N \mhyphen 1}),
	\end{aligned}\\
	&\tilde{\beta}(x;\qn_I) \coloneqq \beta(x;\qn_I)/ \qn_I (\{s\}_0^{N \mhyphen 1}).
\end{align}
$\tilde{\beta}$ is a normalised metric which measures the improvement along the trajectory $\{s\}_0^{N-1}$ relative to the value of the $\qn_I$ there resulting from the addition of a cut at $x$.
In our numerical investigation we observe that using $\beta$ leads to more cuts at points ``far" from the goal states, leading to a richer function approximation there. Using $\tilde{\beta}$, on the other hand, adds cuts closer to the goal states, at points that lead to the highest improvement relative to the value of the $\qn$-function. 
%
%Naturally, there is a trade-off between the richness of the approximation and the cost of extracting a greedy policy online as the computational effort increases with the number of cuts.
% \vspace{-0.2cm}%
\subsection{Advantages of the new formulation}
An efficient implementation that computes the policy \eqref{eq:NstepPi} using an approximation of form \eqref{eq:Qpw_max} uses as its objective function the common terms from all the cuts that build the approximation $\qn_I$ and an epigraph variable $\alpha$ lower bounded by the terms remaining in the cuts. The formulation of the cuts presented in \cite[Lemma 3.1]{menta2020} leads to a MIQCQP due to the extra quadratic terms present in $q_0$, whereas, our formulation \eqref{eq:newlb_defn} gives rise to the following simpler MIQP,
\begin{align*}
    \min_{\{s\}_0^{N \mhyphen 1}} &\sum_{t=0}^{N \mhyphen 1} \left(\begin{aligned}
    \tfrac{\gamma^t}{2}(x_t - x_g)^\top Q& (x_t - x_g) + \\ &\tfrac{\gamma^t}{2}(u_t - u_g)^\top R (u_t - u_g)
    \end{aligned} \right)  + \alpha\\
    \text{s.t.}\quad & \{s\}_0^{N \mhyphen 1} \in \mathcal{S}^{(N)}(x_0)\\
    & \nu_i^\star \left(A x_{N\mhyphen1} + B_1 u_{N\mhyphen1} + B_2 \delta_{N\mhyphen1} + B_3 z_{N\mhyphen1}\right) + c_i \leq \alpha,\\& \qquad \qquad \qquad \qquad \qquad \qquad \qquad \qquad \quad i = 1,\dots,I,
\end{align*}
where the $\nu_i^\star$ is the optimal dual multiplier used in $q_i$. The removal of the quadratic constraints arising from $q_0$ leads to a significant reduction of the policy computation times.

The cuts we generate add linear terms in addition to the  quadratic terms present in cut $q_0$, removing the advantage $q_0$ enjoyed in \cite{menta2020} where the extra quadratic terms dominated the added linear terms, especially, closer to the goal states. In combination with the proposed metric $\tilde{\beta}$, compared to the formulation in \cite{menta2020}, our formulation is able to add more cuts that build the $\qn_I$ closer to the goal states giving a richer approximation there, this leads to better closed-loop performance as the trajectories from any initial condition definitely traverse this region to reach $x_g$.

%%%%%%%%%%%%%%%%%%%%%%%%%%%%%%%%%%%%%%%%%%%%%%%%%%%%%%%%%%%%%%%%%%%%%%%%%%%%%%%%
\section{Numerical Examples}\label{sec:NumEx}
We test the performance of the resulting controller on a traction control problem \cite{traction_control} and a boiler-turbine system \cite{boiler_turbine}.
\subsection{Traction control}
A traction controller is used to improve the stability and steerability of a vehicle under slippery conditions. The controller aims to control the \textit{slip} between the tires and the road to maximise the frictional torque. The frictional torque depends on slip through a non-linear relation, which is approximated as a piecewise affine function with 2 regions, leading to hybrid vehicle dynamics with 2 modes \cite{traction_control}.
\subsubsection{System dynamics}
The dynamics of the vehicle \cite[(7)]{traction_control} are modelled as an MLD system with the following matrices:
\begin{align}\label{eq:traction_control_dyn}
A &=\begin{bmatrix}
1 & 0 & 0 \\
4.8792 \cdot 10^{\mhyphen2} & 1.0005 & \mhyphen2.1835 \cdot 10^{\mhyphen2}\\
\mhyphen1.5695 \cdot 10^{\mhyphen7} & \mhyphen6.4359 \cdot 10^{\mhyphen6} & 1.0003
\end{bmatrix},\\
B_1 &= \begin{bmatrix}
1 & 0\\
4.8792 \cdot 10^{\mhyphen2} & \mhyphen6.5287\\
\mhyphen1.5695 \cdot 10^{\mhyphen7} & 8.9695 \cdot 10^{\mhyphen2}
\end{bmatrix},\\
B_2 &= \begin{bmatrix}
0 & 0\\
0.10943 & 0.81687\\
\mhyphen1.5034 \cdot 10^{\mhyphen3} & \mhyphen1.1223 \cdot 10^{\mhyphen2}
\end{bmatrix},\quad
B_3 = \begin{bmatrix}
0 & 0\\
1 & 0\\
0 & 1
\end{bmatrix}.
\end{align}
The states are $ x = \begin{bmatrix} \tau_c^{\mhyphen} & \omega_e & v_v \end{bmatrix}^\top$, where $\tau_c^{\mhyphen}$ is the combustion torque applied at the previous time step, $\omega_e$ is the engine speed, $v_v$ is the vehicle velocity. The inputs are $ u = \begin{bmatrix} \Delta \tau_c & \mu_s\end{bmatrix}^\top$, where $\Delta \tau_c$ is the change in the combustion torque and $\mu_s$ is the coefficient of friction between the tire and the road, an uncontrolled input that is estimated and given to the controller \cite[Sec. VI.C]{traction_control}. The auxiliary binary variables $\delta = \begin{bmatrix} \delta_1 & \delta_2 \end{bmatrix}^\top$, describe the modes of operation of the vehicle dynamics and the auxiliary continuous variables $z = \begin{bmatrix} z_1 & z_2 \end{bmatrix}^\top$ capture the mode dependent dynamics of the system.
Slip $(\Delta \omega)$ at a given time step is computed as:
\[
\Delta \omega = \frac{v_v}{r_t} - \frac{\omega_e}{g_r} = \underbrace{\begin{bmatrix} 0 & \mhyphen\frac{1}{g_r} & \frac{1}{r_t} \end{bmatrix}}_{S_x} x = S_x x,
\]
where $r_t$ is the tire radius and $g_r$ is the driveline gear ratio. The following constraints apply to the dynamics:
\begin{equation}\label{eq:TCBoxCons}
	\mhyphen20 \leq \tau_c \leq 176, \quad \mhyphen40 \leq \Delta \tau_c \leq 40, \quad	\Delta \omega \geq 0.
\end{equation}
The constraints ensuring well-posedness and \eqref{eq:TCBoxCons} form the MLD constraints for this system where $n_c = 25$. The MLD constraint matrices are given in the Appendix.
\subsubsection{Stage cost function}
It is observed that the frictional torque is the highest at the line separating the 2 regions of operation. We make the assumption that $\mu_s$ is fixed, resulting in a piecewise affine relation between slip and frictional torque. Hence, the control problem is to steer the system to the switching boundary, as that maximises the frictional torque, with a penalty on changing combustion torque. We choose the following quadratic stage cost function:
\begin{align}\label{eq:StageCost}
	&\ell(x, u) = \frac{1}{2} \cdot 50 \cdot (\Delta \omega - \Delta \omega_g)^2 + \frac{1}{2} \cdot 0.1 \cdot (\Delta \tau_c)^2 \nonumber\\
			   &= \frac{1}{2} (x - x_g)^\top Q (x - x_g) + \frac{1}{2}  (u - u_g)^\top R (u - u_g)\\
  	&Q = \begin{bmatrix} 0 & 0 & 0\\ 0 & 0.259 & \mhyphen 12.08\\ 0 & \mhyphen 12.08 & 563.04 \end{bmatrix} \qquad R = \begin{bmatrix} 0.1 & 0\\ 0 & 0 \end{bmatrix}\nonumber
\end{align}
Since we have $ Q\succeq 0$ and $R\succeq 0$, we shift their eigenvalues by a small $\epsilon = 10^{\mhyphen9}$ to make them positive definite, this shift has an negligible impact on performance. We choose $\mu_s = 0.1934$ which results in $\Delta \omega_g = 2$ and any $x_g$ that satisfies $\Delta \omega_g = S_x x_g$ and $u_g = \begin{bmatrix} 0 & \mu_s \end{bmatrix}^\top$.

\begin{table*}[t]
    \centering
    \caption{Comparison of closed-loop costs against HMPC for the traction control example}
    \begin{tabular}{|c|c|c|c|c|c|c|}
        \hline
        \multirow{7}{*}{Baseline}&\multicolumn{3}{c|}{\multirow{2}{*}{$\pi_1:$ HMPC $N=30$}}&\multicolumn{3}{c|}{\multirow{2}{*}{$\pi_1:$ HMPC $N=15$}}\\
        &\multicolumn{3}{c|}{}&\multicolumn{3}{c|}{}\\\cline{2-7}
         & \multirow{5}{1.6cm}{$\%$ of points where $\pi_1$ is better than the baseline} & \multirow{5}{2cm}{Mean better performance $(\%)$ at points where $\pi_1$ is better} & \multirow{5}{2cm}{Mean worse performance $(\%)$ at points where $\pi_1$ is worse} & \multirow{5}{1.6cm}{$\%$ of points where $\pi_1$ is better than the baseline} & \multirow{5}{2cm}{Mean better performance $(\%)$ at points where $\pi_1$ is better} & \multirow{5}{2cm}{Mean worse performance $(\%)$ at points where $\pi_1$ is worse}\\
         & & & & & &\\
         & & & & & &\\
         & & & & & &\\
         & & & & & &\\\hline
        $\pi(x;\mathcal{Q}^{(15)}_{22})$ & 48.7 & 0.12 & 0.08 & 25.7 & 0.11 & 1.53 \\[2.5pt]
        HMPC $N=15$                          & 90.2 & 1.26 & 0.13 & - & - & -\\\hline
    \end{tabular}
    \label{tab:HMPC_TC}
\end{table*}

\subsubsection{$\mathcal{Q}^{(N)}$-function approximation}
The $\qn$-function is learnt with $N=15$, over $\xalgm$ with $22$ data points spread over the slip range $[0,15]$, $I_{max} = 100$ and $\eta_{min} = 10^{\mhyphen 5}$. The generated approximation $\mathcal{Q}^{(15)}_{22}$ has $I = 22$.

\begin{figure}[t]
% 	\vspace{-0.3cm}
	\centering
	\title{Traction control system}
	\begin{minipage}{.55\columnwidth}
		\centering
		\includegraphics[width=\textwidth]{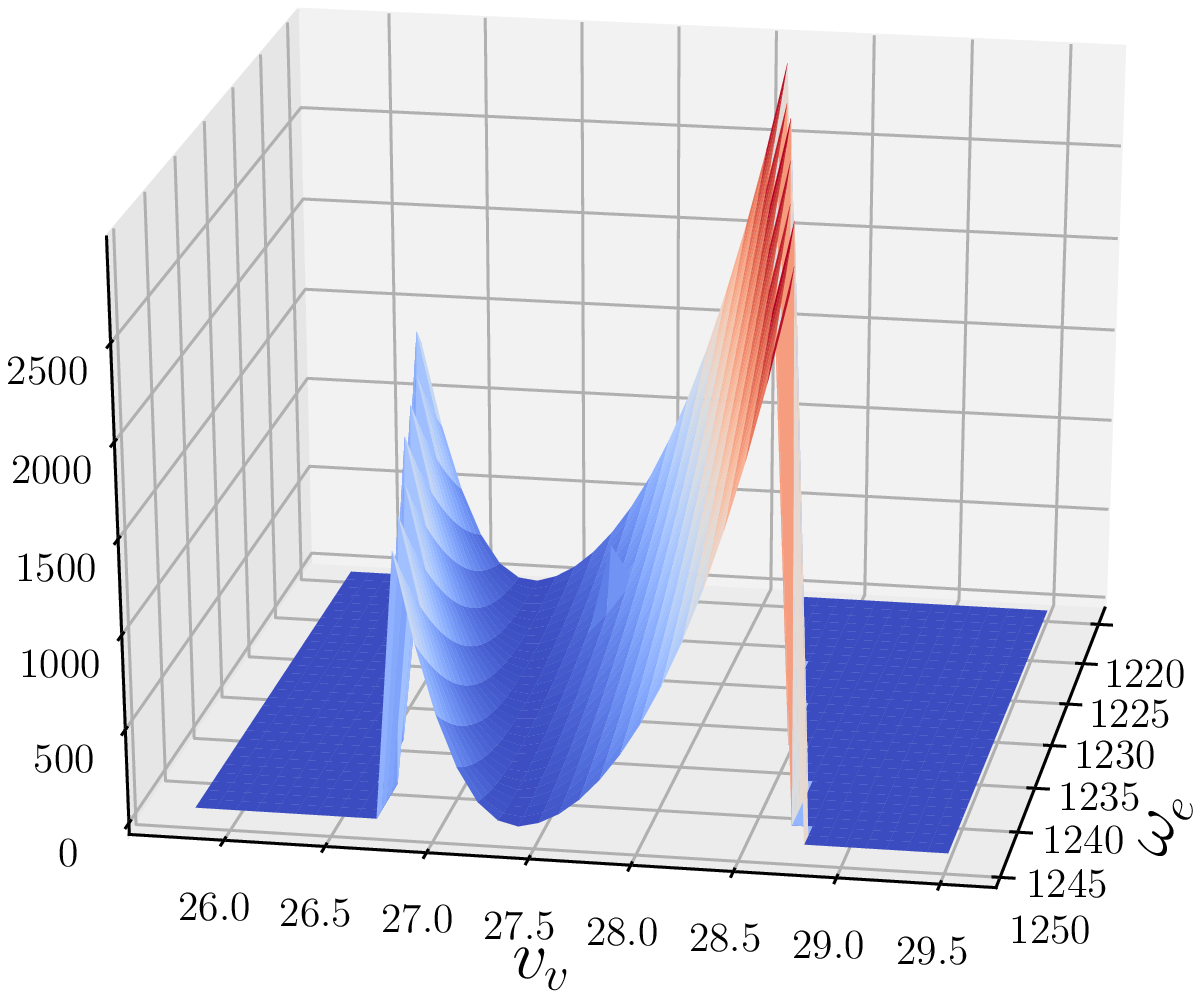}
		\caption{Plot of $V^\star(x)$ for the \textcolor{white}{lotslotslots lotslotsl} traction control problem}
		\label{fig:TCVPlot}
	\end{minipage}%
	\begin{minipage}{.55\columnwidth}
		\hspace{-0.9cm}\centering
		\includegraphics[width=\textwidth]{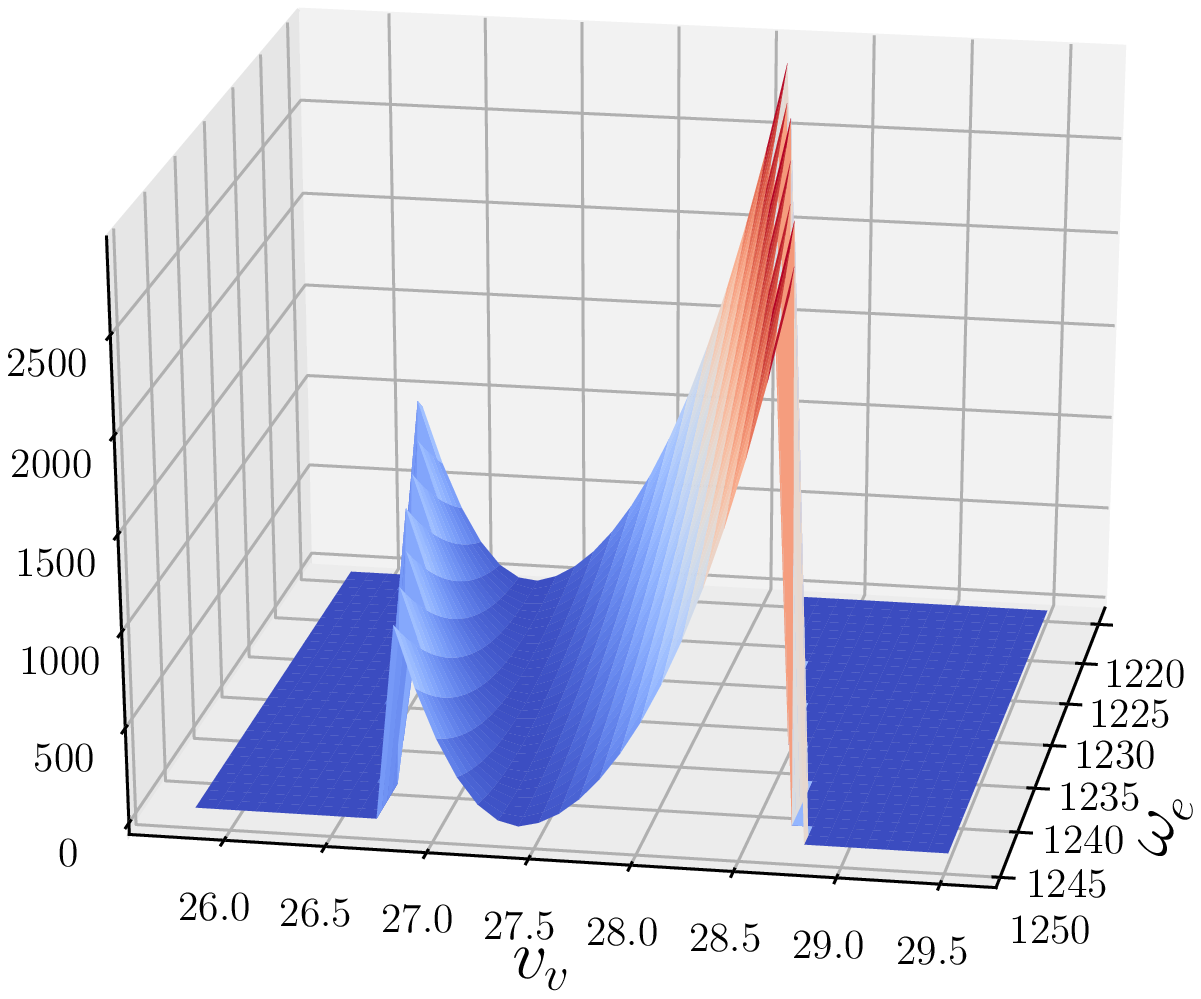}
		\caption{Plot of $\underbar{\text{$\mathcal{Q}$}}^{(15)}_{\,22}(x)$ for the \textcolor{white}{lotslots} traction control problem}
		\label{fig:TCQPlot}
	\end{minipage}
% 	\vspace{-0.4cm}
\end{figure}

As the system dimension is small, we approximate $V^\star(x)$ by using HMPC with a long horizon length of $30$ for the points in the slip range $[0,6.5]$. For the purposes of visualisation of the trained $\qn$-functions we define
\[
\textstyle \barqn(x) = \inf_{\{s\}_0^{N-1} \in \mathcal{S}^{(N)}(x)} \qn(\{s\}_0^{N-1}).
\] Evaluating $\barqn$ gives us an idea of the approximation quality as ${\barqn}^\star(x) = V^\star(x)$. The plots of $V^\star$ and the $\underbar{\text{$\mathcal{Q}$}}^{(15)}_{\,22}$ are shown in Figs. \ref{fig:TCVPlot} and \ref{fig:TCQPlot}, for the slip range $[0,6.5]$. Figs. \ref{fig:TCVPlot} and \ref{fig:TCQPlot} suggest that the $\qn$-function approximates the optimal value function very well.

\subsubsection{Simulations using $\mathcal{Q}^{(N)}_I$}
We compare the closed-loop costs of applying the control policy \eqref{eq:NstepPi}, derived from $\mathcal{Q}^{(15)}_{22}$, against those of HMPC, with the LQR terminal cost derived using $A$ and $B_1$ matrices, from $470$ initial conditions spread across the slip range $[0,6.5]$. 

\begin{figure}[b]
% 	\vspace{-0.3cm}
	\centering
	\title{Traction control system}	\includegraphics[width=\columnwidth]{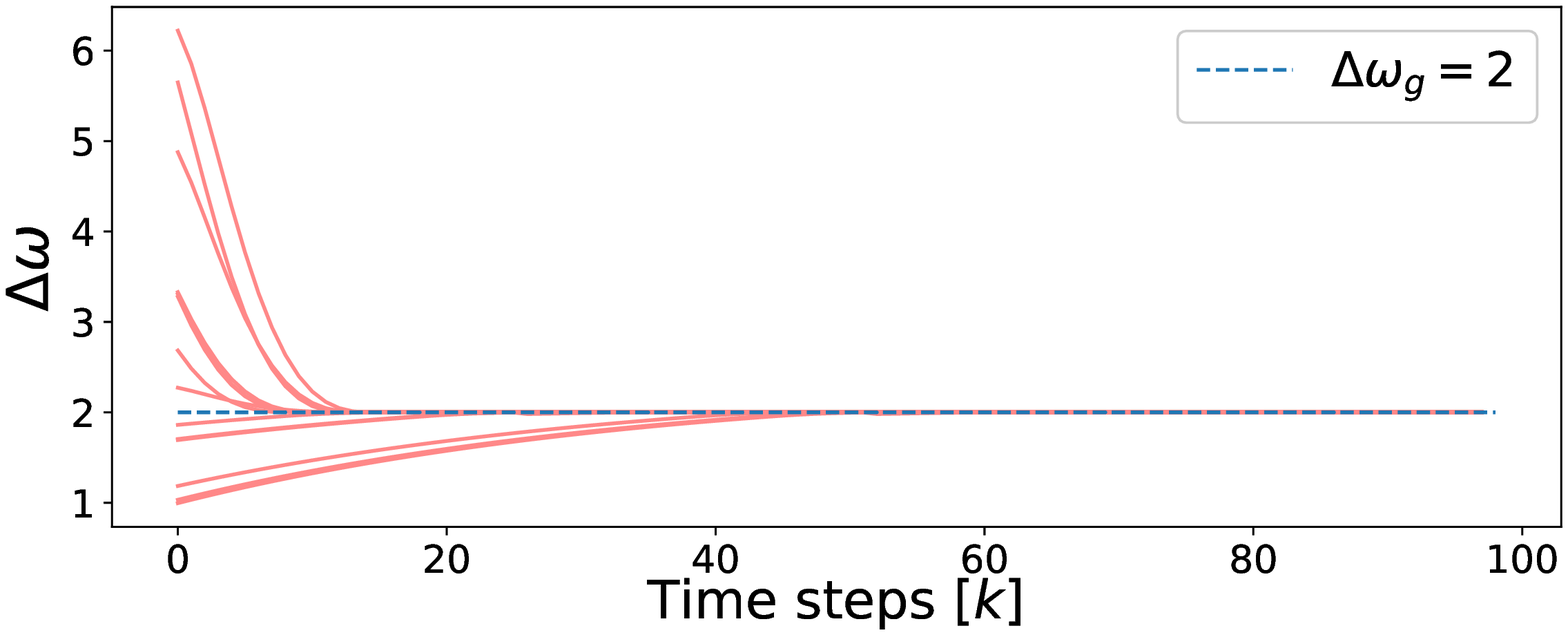}
    \caption{Slip trajectories generated by $\pi(x;\mathcal{Q}^{(15)}_{22})$}
    \label{fig:TC_MPC_15}
    \vspace{-0.2cm}
\end{figure}

\textcolor{black}{Table \ref{tab:HMPC_TC} presents a detailed comparison of the closed-loop costs induced by the different controllers from the $470$ initial conditions.} The closed-loop costs of our method and HMPC $N = 15$ are on an average $0.02\%$ and $1.1\%$ worse than the values of $V^\star(x)$ respectively, which shows that our control policy is effectively identical to the optimal control policy. \textcolor{black}{ Our controller $\pi(x;\mathcal{Q}^{(15)}_{22})$ clearly performs better than HMPC of the same horizon length, as seen from the entries of Table \ref{tab:HMPC_TC}.} \textcolor{black}{The points where HMPC $N=15$ and our controller outperform HMPC $N=30$, it is observed that they do so by very small values due to numerical effects. The $\Delta \omega$ trajectories resulting from $\pi(x;\mathcal{Q}^{(15)}_{22})$ are shown in Fig. \ref{fig:TC_MPC_15}.} The average time taken to compute the policy across the $470$ initial conditions is $~0.132s$ for HMPC $N=15$, $~0.141s$ for $\pi(x;\mathcal{Q}^{(15)}_{22})$ and $0.394s$ for HMPC $N=30$. For the sake of comparison, solving the MIQCQP to compute the controller of \cite{menta2020} for this system using an approximation with the same $N$ and $I$ would require in excess of $0.25s$. On computing $\mathcal{Q}^{(6)}_{22}$ at $2327$ points in the slip range $[0,6.5]$, it is seen that the $q_0$ is used at only $34$ points, whereas, an approximation trained with the same parameters and metric using the formulation in \cite{menta2020} uses $q_0$ at $1357$ of these points.

We also tested the effect of using $\pi(x; \mathcal{Q}^{(15)}_{22})$ on the system when $\mu_s$ was different than the value used in training. \textcolor{black}{We observed that when the difference between $\Delta \omega_g$ and \textit{slip at the switching boundary} is not large our controller drives the slip of the system to the one closer to the initial state.} Since $\mu_s$ can change dramatically, one can train $\mathcal{Q}$-functions for various values of $\mu_s$ and control the vehicle using the $\mathcal{Q}$-function trained with the $\mu_s$ closest to the estimated one.

\begin{figure}[t]
% 	\vspace{-0.3cm}
	\centering
	\title{Traction control system}
	\includegraphics[width=\columnwidth]{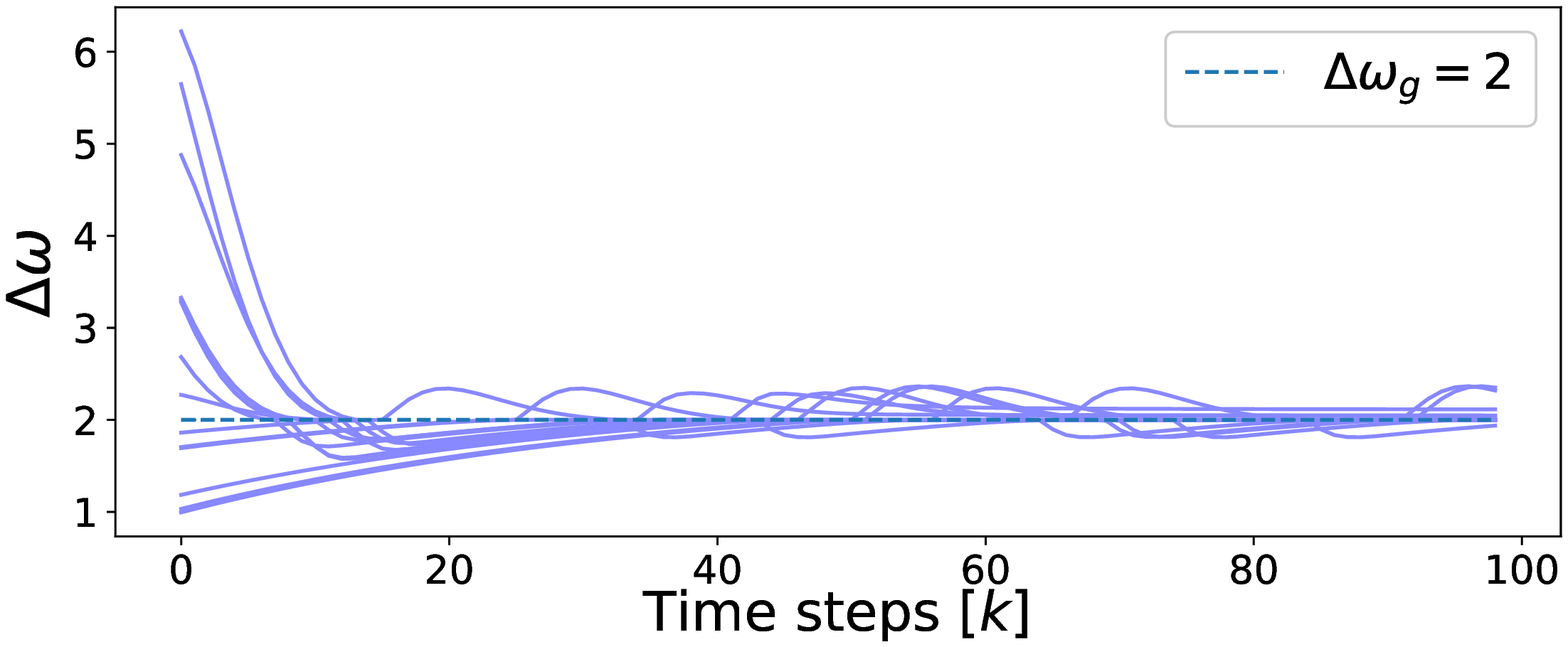}
    \caption{Slip trajectories generated by the best NN model}
    \label{fig:TC_NN_15}
    \vspace{-0.5cm}
\end{figure}

To compare against model-free techniques we train a Fully Connected Neural Network with a Sigmoid activation at the last layer and ReLU activations elsewhere using data from HMPC $N=30$. The network predicts $\tau_c$ using $\omega_e, v_v, \Delta \omega, \tau_c^\mhyphen$ and the region number. We used BOHB \cite{BOHB}
%with a max. budget for the number of iterations of $2000$, min. budget of $50$ and $\eta = e$, 
to find a good architecture searching over models with $4$-$7$ layers with $8$-$64$ neurons each and over the range $[10^{\mhyphen 1}, 10^{\mhyphen 4}]$ for the learning rate. The models were trained and tested using $75927$ and $31637$ data points respectively. The best model had $2,776$ parameters and gave a mean test error of $0.186\%$ with a variance of $0.223\%$. The resulting controller, resulted in large oscillations at the switching boundary due to the high sensitivity of the dynamics to $\tau_c$ there. \textcolor{black}{Fig. \ref{fig:TC_NN_15} shows the trajectories resulting from the use of the NN controller from the same initial states as in Fig. \ref{fig:TC_MPC_15}.} This suggests that off-the shelf neural network techniques might not be the best choice for hybrid problems due to the switching dynamics compared to other methods that approximate explicit MPC using model information \cite{Chen_neural}, \cite{maddalena2019neural}.

\subsection{Boiler-turbine system}
A boiler-turbine system converts chemical energy from burning fuels to electrical energy via a steam boiler and a turbine. The boiler-turbine system has been widely studied with a variety of control techniques. This system has several nominal operating points \cite{boiler_turbine} and we want to steer the system to the ``normal" nominal operating point starting from some other state.

\subsubsection{System dynamics}
The non-linear dynamics of the system \cite[eqn. (3), (5)]{boiler_turbine} are given by
\begin{subequations}
	\begin{align}
	\frac{dp}{dt} &= \mhyphen 0.0018 u_2 p^{9/8} + 0.9u_1 - 0.15 u_3,\\
	\frac{dp_0}{dt} &= (0.073 u_2 - 0.016)p^{9/8} - 0.1p_0,\\
	\frac{dp_f}{dt} &= \frac{141 u_3 - (1.1 u_2 - 0.19)p}{85},\\
	X_w &= 0.05\left(0.1307 p_f + 100 \alpha_{cs} + \frac{q_e}{9} - 67.975\right),\\
	\alpha_{cs} &= \frac{(1 - 0.001538p_f)(0.8p - 25.6)}{p_f (1.0394 - 0.0012304p)},\\
	q_e &= (0.854 u_2 - 0.147)p + 45.59u_1 \nonumber\\
	& \qquad \qquad \qquad \qquad \qquad - 2.514 u_3 - 2.096,
	\end{align}
\end{subequations}
	where $p$ is the drum pressure (kg/cm$^2$), $p_0$ the power output (MW) and $p_f$ the fluid density (kg/cm$^3$). The inputs $u_1$, $u_2$ and $u_3$ respectively control fuel flow valve position, steam control valve position and feed-water flow valve position and lie in the interval $[0,1]$. $X_w$, $q_e$ and $\alpha_{cs}$ denote the drum water level (m), the evaporation rate (kg/s) and the steam quality, respectively. Conversion of the non-linear dynamics to an MLD system is done as described in Section 2 of \cite{boiler_turbine}.
	%The non-linear terms $u_2 p^{9/8}, u_2 p, p^{9/8}$ and $\alpha_{cs}$ are linearised by splitting the operation ranges of $p$, $p_f$ and $u_2$ into 4 equal sections each and then fitting linear functions to them using the least squares fitting technique in each of the operating regions \cite{boiler_turbine}. The PWA dynamics that result from the linearisation are converted to a MLD system by using HYSDEL \cite{HYSDEL} and then discretised with a sampling interval of $T_s$.
	The resulting MLD system has $7$ states - $p[k], p_0[k], p_f[k], X_w[k\mhyphen1], u_1[k\mhyphen1], u_2[k\mhyphen1]$ and $u_3[k\mhyphen1]$, $3$ inputs, $64$ continuous and $45$ binary auxiliary variables and $n_c=414$. The system is subject to the following constraints:
	\begin{multline*}
		\left|\frac{u_1[k] - u_1[k\mhyphen1]}{T_s}\right| \leq 0.007 , \quad \left|\frac{u_3[k] - u_3[k\mhyphen1]}{T_s}\right| \leq 0.05,\\ \mhyphen2 \leq \frac{u_2[k] - u_2[k\mhyphen1]}{T_s} \leq 0.02,  \left|\frac{X_w[k] - X_w[k\mhyphen1]|}{T_s}\right| \leq 0.02.
	\end{multline*}

\subsubsection{Stage cost function} The stage cost depends only on the state with $Q = diag([742,441,1,\epsilon,\epsilon,\epsilon,\epsilon])$ where $\epsilon = 10^{\mhyphen 9}$ and $R$ a $3 \times 3$ zero matrix. The operations in computing the cuts are modified to remove the dependency on $R$. $x_g$ is the \textit{normal} operating point $p = 108$, $p_0 = 66.65$, $p_f = 428$ along with the corresponding equilibrium inputs $u_1 = 0.34$, $u_2 = 0.684$, $u_3 = 0.435$ and $X_w = 0.0129$.

\subsubsection{$\qn$-function approximation} Let us denote the operational ranges of $p,p_0$ and $p_f$ using $\mathcal{R}_p, \mathcal{R}_{p_0}$ and $\mathcal{R}_{p_f}$. The size of this system makes it computationally expensive to train $\qn$-functions for long horizon lengths. We run Alg. \ref{alg:QBenders_xudz} using $N=6$ over $\xalgm$ with $256$ data points distributed uniformly in $\mathcal{R}_{p}\times\mathcal{R}_{p_0}\times\mathcal{R}_{p_f}$, using $I_{max} = 200$ and $\eta_{min} = 10^{\mhyphen3}$. The trained approximation, $\mathcal{Q}^{(6)}_{200,\tilde{\beta}}$, has $I=200$. \textcolor{black}{We build the approximations $\mathcal{Q}^{(6)}_{I,\tilde{\beta}}$ using the first $I = 25, 50, 75, 100, 125$ and $150$ cuts.} To demonstrate the advantage of the improvement metric $\tilde{\beta}$ over $\beta$ we train an approximation by running Alg. \ref{alg:QBenders_xudz} using the improvement metric $\beta$ instead of $\tilde{\beta}$. Setting $I_{max} = 150$ and $\eta_{min} = 10^{\mhyphen3}$, we train the approximation $\mathcal{Q}^{(6)}_{150,\beta}$ with $I=150$ cuts. \textcolor{black}{We build the approximations $\mathcal{Q}^{(6)}_{I,\beta}$ using the first $I = 25, 50, 75, 100$ and $125$ cuts.}
% For this system, $\mathcal{Q}^{(6)}_{100, B}$ is richer closer to $x_g$ and $\mathcal{Q}^{(6)}_{100, A}$ is richer at points far from $x_g$.
% In an attempt ot make the best of both approximations we combine the best cuts from the two algorithms to build the approximation $\mathcal{Q}^{(6)}_{85, C}$, this is done by evaluating $\mathcal{Q}^{(6)}_{100, A}$ and $\mathcal{Q}^{(6)}_{100, B}$ at $2592$ points distributed uniformly across $\mathcal{R}_{p}\times\mathcal{R}_{p_0}\times\mathcal{R}_{p_f}$ and only retaining cuts which are active in at least $5$ points.

\subsubsection{Simulations using $\mathcal{Q}^{(N)}_I$}
We compare the closed-loop costs of HMPC, with the LQR terminal cost which just uses the $3\times3$ principal submatrix of $Q$, against those of applying \eqref{eq:NstepPi} derived from \textcolor{black}{$\mathcal{Q}^{(6)}_{I,\tilde{\beta}}$ and $\mathcal{Q}^{(6)}_{I,\beta}$, for various values of $I$} from $100$ initial conditions spread uniformly across $\mathcal{R}_{p}\times\mathcal{R}_{p_0}\times\mathcal{R}_{p_f}$ in Table \ref{tab:HMPC}. The performance of our controllers is close to that of HMPC $ N = 6$ with our controllers outperforming it at several points. The average time taken to compute the policy is $0.6s$ for HMPC $N=6$, $0.88s$ for $\pi(x;\mathcal{Q}^{(6)}_{100,\beta})$, $0.88s$ for $\pi(x;\mathcal{Q}^{(6)}_{100,\tilde{\beta}})$, $1.01s$ for $\pi(x;\mathcal{Q}^{(6)}_{150,\tilde{\beta}})$, $1.23s$ for $\pi(x;\mathcal{Q}^{(6)}_{200,\tilde{\beta}})$ and $21.78s$ for HMPC $N=20$. The controller derived from our approximations and HMPC $N = 6$ have closed-loop costs within $2.5\%$ of HMPC $N = 20$ at a fraction of the computation cost.

\textcolor{black}{For $I = 25, 50, 75$ and $100$ the closed-loop performance of $\pi(x;\mathcal{Q}^{(6)}_{I,\tilde{\beta}})$ was seen to be better than that of $\pi(x;\mathcal{Q}^{(6)}_{I,\beta})$ but the performance was similar for $I = 125$ and $150$, so using approximations built with the metric $\tilde{\beta}$ is better for small $I$.} For $N=6$ and $I=100$, solving the MIQCQP controller of \cite{menta2020} for this system would require in excess of $1.25s$, $42\%$ higher than our controller. On computing $\mathcal{Q}^{(6)}_{100,\tilde{\beta}}$ at $2592$ points distributed uniformly across $\mathcal{R}_{p}\times\mathcal{R}_{p_0}\times\mathcal{R}_{p_f}$, it is seen that $q_0$ is not at all used, whereas, an approximation trained with the same parameters and metric using the formulation in \cite{menta2020} uses $q_0$ at $1812$ of these points.

\begin{table*}[t]
    \centering
    \caption{Comparison of closed-loop costs against HMPC for the boiler-turbine system ($\qn_I$ is an approximation of horizon length $N$ with $I$ cuts)}
    \begin{tabular}{|c|c|c|c|c|c|}
        \hline
        \multirow{6}{*}{Baseline}&\multicolumn{2}{c|}{\multirow{2}{*}{$\pi_1:$ HMPC $N=20$}}&\multicolumn{3}{c|}{\multirow{2}{*}{$\pi_1:$ HMPC $N=6$}}\\
        &\multicolumn{2}{c|}{}&\multicolumn{3}{c|}{}\\\cline{2-6}
         & \multirow{4}{2.1cm}{$\%$ of points where $\pi_1$ is better than the baseline} & \multirow{4}{2.4cm}{Mean better perfor- mance $(\%)$ at points where $\pi_1$ is better} & \multirow{4}{2.1cm}{$\%$ of points where $\pi_1$ is better than the baseline} & \multirow{4}{2.4cm}{Mean better perfor- mance $(\%)$ at points where $\pi_1$ is better} & \multirow{4}{2.4cm}{Mean worse perfor- mance $(\%)$ at points where $\pi_1$ is worse}\\
         & & & & &\\
         & & & & &\\
         & & & & &\\\hline
        $\pi(x;\mathcal{Q}^{(6)}_{25,\beta})$ & 100 & 2.93 & 96 & 0.92 & 0.01 \\[2.5pt]
        $\pi(x;\mathcal{Q}^{(6)}_{50,\beta})$ & 100 & 2.73 & 90 & 0.84 & 0.56 \\[2.5pt]
        $\pi(x;\mathcal{Q}^{(6)}_{75,\beta})$ & 100 & 2.41 & 84 & 0.56 & 0.54 \\[2.5pt]
        $\pi(x;\mathcal{Q}^{(6)}_{100,\beta})$ & 100 & 2.46 & 90 & 0.59 & 0.92 \\[2.5pt]
        $\pi(x;\mathcal{Q}^{(6)}_{125,\beta})$ & 100 & 2.24 & 81 & 0.38 & 0.44 \\[2.5pt]
        $\pi(x;\mathcal{Q}^{(6)}_{150,\beta})$ & 100 & 2.09 & 70 & 0.22 & 0.27 \\[2.5pt]
        $\pi(x;\mathcal{Q}^{(6)}_{25,\tilde{\beta}})$ & 100 & 2.65 & 86 & 0.73 & 0.09 \\[2.5pt]
        $\pi(x;\mathcal{Q}^{(6)}_{50,\tilde{\beta}})$ & 100 & 2.6 & 84 & 0.7 & 0.16 \\[2.5pt]
        $\pi(x;\mathcal{Q}^{(6)}_{75,\tilde{\beta}})$ & 100 & 2.33 & 73 & 0.44 & 0.08 \\[2.5pt]
        $\pi(x;\mathcal{Q}^{(6)}_{100,\tilde{\beta}})$ & 100 & 2.28 & 75 & 0.38 & 0.11 \\[2.5pt]
        $\pi(x;\mathcal{Q}^{(6)}_{125,\tilde{\beta}})$ & 100 & 2.31 & 76 & 0.4 & 0.06 \\[2.5pt]
        $\pi(x;\mathcal{Q}^{(6)}_{150,\tilde{\beta}})$ & 100 & 2.16 & 66 & 0.26 & 0.09 \\[2.5pt]
        $\pi(x;\mathcal{Q}^{(6)}_{200,\tilde{\beta}})$ & 100 & 2.07 & 57 & 0.22 & 0.16 \\[2.5pt]
        HMPC $N=6$                          & 100 & 1.92 & - & - & -\\[2.5pt]\hline
    \end{tabular}
    \label{tab:HMPC}
\end{table*}

%%%%%%%%%%%%%%%%%%%%%%%%%%%%%%%%%%%%%%%%%%%%%%%%%%%%%%%%%%%%%%%%%%%%%%%%%%%%%%%%
\section{Conclusions}
We presented a new formulation for the lower bounding functions that are used in the $\qn$-function approximations along with a modified algorithm that improves the performance of our controller. We show that our controller can achieve good control performance on real-world examples through simulations on the traction control problem and the boiler-turbine system. \textcolor{black}{Our algorithm is shown to work well for systems with dimensions relevant for practical MPC problems} with the performance of our controller comparable to that of HMPC with an equal horizon length. The closed-loop costs of our controller also come close to the costs of the optimal control policy (approximated by long horizon length HMPC).% while also learning interesting details of $\oqn$.

%%%%%%%%%%%%%%%%%%%%%%%%%%%%%%%%%%%%%%%%%%%%%%%%%%%%%%%%%%%%%%%%%%%%%%%%%%%%%%%%
\section{Acknowledgements}

Research supported by the European Research Council (ERC) under the European Union’s Horizon 2020 research and innovation programme, grant agreement OCAL, No. 787845.

\bibliographystyle{ieeetr}
\bibliography{refs}
\textcolor{white}{Final line}
\clearpage
% \newpage
\begin{strip}
\appendix{
The MLD constraint matrices corresponding to the traction control system \eqref{eq:traction_control_dyn} are given here:}
\newline
\begin{equation*}
E_1 = \begin{bmatrix}
0 & 5.3781\\
0 & \mhyphen5.3781\\
0 & 0\\
0 & 0\\
0 & 0\\
0 & 0\\
\mhyphen0.000424 & 6.17455\\
0.000424 & \mhyphen6.17455\\
0 & 0\\
0 & 0\\
5.82575 \cdot 10^{\mhyphen6} &  \mhyphen0.0848295\\
\mhyphen5.82575 \cdot 10^{\mhyphen6} & 0.0848295\\
0 & 5.3781\\
0 \mhyphen5.3781\\
\mhyphen1 & 0\\
1 & 0\\
\mhyphen1 & 0\\
0 & 0\\
0 & \mhyphen1\\
0 & 1\\
0 & 0\\
0 & 0\\
0 & 0\\
0 & 0\\
1 & 0
\end{bmatrix},\qquad\qquad\qquad\qquad\qquad\qquad
E_2 = \begin{bmatrix}
0 & \mhyphen73.456664601\\
0 & 57.273443009\\
1 & 1\\
\mhyphen1 \mhyphen1\\
\mhyphen88.586529999 & 0\\
\mhyphen59.558523999 & 0\\
59.558523999 & 0\\
88.586529999 & 0\\
\mhyphen0.817065082 & 0\\
\mhyphen1.216723605 & 0\\
1.216723605 & 0\\
0.817065082 & 0\\
73.456664601 & 0\\
\mhyphen57.273443009 & 0\\
0 & 0\\
0 & 0\\
0 & 0\\
0 & 0\\
0 & 0\\
0 & 0\\
0 & 0\\
0 & 0\\
0 & 0\\
0 & 0\\
0 & 0
\end{bmatrix}
\end{equation*}

\begin{equation*}
E_3 = \begin{bmatrix}
0 & 0\\
0 & 0\\
0 & 0\\
0 & 0\\
1 & 0\\
\mhyphen1 & 0\\
1 & 0\\
\mhyphen1 & 0\\
0 & 1\\
0 & \mhyphen1\\
0 & 1\\
0 & \mhyphen1\\
0 & 0\\
0 & 0\\
0 & 0\\
0 & 0\\
0 & 0\\
0 & 0\\
0 & 0\\
0 & 0\\
0 & 0\\
0 & 0\\
0 & 0\\
0 & 0\\
0 & 0
\end{bmatrix},\qquad
E_4 = \begin{bmatrix}
0 & 0.152993521 & \mhyphen0.713114094\\
0 & \mhyphen0.152993521 & 0.713114094\\
0 & 0 & 0\\
0 & 0 & 0\\
0 & 0 & 0\\
0 & 0 & 0\\
\mhyphen0.000424 & \mhyphen0.173399999 & 0.806695\\
0.000424 & 0.173399999 & \mhyphen0.806695\\
0 & 0 & 0\\
0 & 0 & 0\\
5.82575 \cdot 10^{\mhyphen6} & 0.002377759  & \mhyphen0.011079999\\
\mhyphen5.82575 \cdot 10^{\mhyphen6} & \mhyphen0.002377759 &  0.011079999\\
0 & 0.152993521 & \mhyphen0.713114094\\
0 & \mhyphen0.152993521 & 0.713114094\\
\mhyphen1 & 0 & 0\\
1 & 0 & 0\\
0 & 0 & 0\\
0 & \mhyphen0.719942405 & 3.355704698\\
0 & 0 & 0\\
0 & 0 & 0\\
0 & \mhyphen10 & 0\\
0 & 10 & 0\\
0 & 0 & \mhyphen1\\
0 & 0 & 1\\
0 & 0 & 0
\end{bmatrix},\qquad
E_5 = \begin{bmatrix}
-0.61532\\
57.888762009\\
1\\
-1\\
0\\
0\\
59.558523999\\
88.586529999\\
0\\
0\\
1.216723605\\
0.817065082\\
72.841344601\\
0.615319\\
176\\
20\\
40\\
0\\
0.193439319\\
-0.193439319\\
2500\\
100\\
100\\
20\\
40
\end{bmatrix}
\end{equation*}
\end{strip}

\end{document}